\documentclass[10pt,a4paper]{article}
\usepackage{amsfonts}
\usepackage{amssymb}
\usepackage{mathrsfs}
\usepackage{amsthm}
\usepackage{setspace}
\usepackage{color}
\usepackage{bm}

\usepackage{fancyhdr}
\newtheorem{theo}{Theorem}[section]
\newtheorem{prop}[theo]{Proposition}
\newtheorem{cor}[theo]{Corollary}
\newtheorem{lemm}[theo]{Lemma}

\newtheorem{con}{Conjecture}

\theoremstyle{definition}
\newtheorem*{rem}{Remark}
\newtheorem{defi}[theo]{Definition}

\newenvironment{lproof}{\emph{Proof of Lemma.}}{ \qed \par}

\newcommand{\be}{\begin{eqnarray*}}
\newcommand{\ee}{\end{eqnarray*}}
\newcommand{\beqa}{\begin{eqnarray}}
\newcommand{\eeqa}{\end{eqnarray}}
\newcommand{\ba}{\begin{array}}
\newcommand{\ea}{\end{array}}
\newcommand{\onab}{\overrightarrow{\nabla}}
\newcommand{\mc}{\mathcal}
\newcommand{\mf}{\mathfrak}

\newcommand{\rP}{\mathsf{P}}

\newcommand{\mbb}{\mathbb}
\newcommand{\wt}{\widetilde}

\newcommand{\A}{\mc{A}}

\begin{document}

\title{Non-regular $|2|$-graded geometries I: general theory}
\author{Stuart Armstrong, \\ St Cross College, \\ Oxford, OX1 3LZ, UK}
\date{2009}
\maketitle

\begin{abstract}
This paper analyses non-regular $|2|$-graded geometries, and show that they share many of the properties of regular geometries -- the existence of a unique normal Cartan connection encoding the structure, the harmonic curvature as obstruction to flatness of the geometry, the existence of the first two BGG splitting operators and of (in most cases) invariant prolongations for the standard Tractor bundle $\mc{T}$. Finally, it investigates whether these geometries are determined entirely by the distribution $H = T_{-1}$ and concludes that this is generically the case, up to a finite choice, whenever $H^1(\mf{g}^1,\mf{g})$ vanishes in non-negative homogeneity.
\end{abstract}

\section{Introduction}
For a group $G$ with subgroup $P$, $G$ is naturally a $P$-principal bundle over the homogeneous space $G/P$, and carries the canonical Mauer-Cartan form, a one-form $\omega$ on $G$ with values in the Lie algebra $\mf{g}$ of $G$, given by left translations and by the identification $TG_{Id} = \mf{g}$.

The curved analogue of such a connection involves a manifold $M$ of same dimension as $G/P$ with a principal $P$-bundle $\mc{P}$ and a one-form $\omega$ which is a section of $T\mc{P}^* \otimes \mf{g}$, subject to certain properties.

The groups $(G,P)$ form a parabolic pair if $G$ is semisimple and $P$ is a parabolic subgroup of $G$ -- i.e.~if the Lie algebra $\mf{g}$ of $G$ is $|k|$-graded and the Lie algebra $\mf{p}$ of $P$ consists of those elements of $\mf{g}$ of non-negative grading.

Cartan connections on for parabolic pairs $(G,P)$ have many elegant properties (\cite{TCPG}, \cite{TBIPG}) -- such as the existence of a filtration of the tangent bundle $T$ of $M$:
\be
T = T^{-k} \supset T^{-k+1} \supset \ldots \supset T^{-2} \supset T^{-1},
\ee
and an algebra\"ic bracket $\mc{K}$ on the associated graded bundle
\be
gr(T) = (T^{-k}/T^{-k+1}) \oplus (T^{-k+1}/T^{-k+2}) \oplus \ldots \oplus (T^{-2}/T^{-1}) \oplus T^{-1}.
\ee
If this filtration has the further property that $[T^j,T^i] \subset [T^{j+i}]$, then there is also another natural bracket on the graded bundle $gr(T)$, the Levi bracket $\mc{L}$, derived from the Lie bracket.

If $\mc{L} = \mc{K}$, then the geometry is said to be regular. Such regular parabolic geometries have an incredibly elegant theory, see \cite{old1926}, \cite{old1931}, \cite{ein}, \cite{old1947}, \cite{tanaka} and \cite{old1994}, for a historical overview, and \cite{CapWeyl}, \cite{CartEquiv}, \cite{paradef}, \cite{TCPG}, \cite{TBIPG} and \cite{capslo} for a summary of modern results. The book \cite{capslo} is particularly good for summarising the constructions and results.

The curvature of a Cartan connection $\omega$ is a section $\kappa$ of $\wedge^2 T^* \otimes \A$ for $\A$ the so-called adjoint bundle of $\mc{P}$. There is a natural cohomology operator $\partial^* : \wedge^2 T^* \otimes \A \to T^* \otimes \A$. Then a Cartan connection $\omega$ is called normal if
\be
\partial^* \kappa = 0.
\ee
The $|k|$-grading of $\mf{g}$ allows us to define the subalgebra $\mf{g}^1 \subset \mf{p} \subset \mf{g}$ consisting of those elements of $\mf{g}$ of strictly positive grading. This leads to the cohomology spaces $H^1(\mf{g}^1,\mf{g})$, which have a splitting according to homogeneity (derived from the grading of $\mf{g}$).

The most important result is that if $H^1(\mf{g}^1,\mf{g})$ vanishes in strictly positive homogeneities (a highly generic condition), there is a unique regular normal Cartan connection compatible with the underlying geometry.

But what is the underlying geometry precisely? This can also be answered by the cohomology spaces, the result being that if $H^1(\mf{g}^1,\mf{g})$ vanishes in non-negative homogeneities (a highly generic condition for $|k|$-grading with $k \geq 2$), the underlying geometry is derived entirely from the regular filtration of the tangent bundle.

How does this extend to non-regular geometries? This paper is dedicated to exploring non-regular geometries in the simplest setting where they can occur: that of $|2|$-graded parabolic geometries. In this grading, the filtration is simply
\be
T = T^{-2} \supset T^{-1},
\ee
so the only data is the distribution $T^{-1}$, which will be designated by $H$.

A non-regular $|2|$-graded geometry is given by a manifold $M$ with tangent bundle $T$ and a distribution $H \subset T$, a $|2|$-graded parabolic pair $(G,P)$ with compatible dimensions, and a principal $P$-bundle $\mc{P} \to M$ that is a structure bundle for $gr(T) = H \oplus T/H$. This reduction in structure group gives the algebra\"ic bracket $\mc{K}$ on $gr(T)$ and defines $\partial^*$ in homogeneity zero. Then the geometry is said to be partially regular if
\be
\partial^* (\mc{K} - \mc{L}) = 0.
\ee
Then there is a similar result for Cartan connections in this setting:
\begin{theo}
Given a partially regular $|2|$-graded geometry $(M,\mc{P})$ such that $H^1(\mf{g}^1,\mf{g})$ vanishes in strictly positive homogeneities, there exists a unique normal Cartan connection compatible with the structure.
\end{theo}

But what if the original data is not given by $H$ and $\mc{P}$, but simply by $H$? We call such an $H$ partially regularisable for a parabolic pair $(G,P)$ if there exists a $P$-principal bundle $\mc{P} \to M$ making this into a partially regular geometry.

If $H$ is holomorphic and generic, then the result is that:
\begin{theo}
A generic holomorphic distribution $H$ of correct rank and co-rank for a $|2|$-graded parabolic pair $(G,P)$ is partially regularisable almost everywhere. If $H^1(\mf{g}^1,\mf{g})$ vanishes in cohomology zero, there is only a finite choice of such partial regularisations. There exists uniqueness procedures that allow one to pick a single partial regularisation in a well defined fashion.
\end{theo}
Similar results apply for real distributions $H$ that are sufficiently close to a regular distribution. Whether the geometry is partially regularisable at a point depends only on the Levi bracket $\mc{L}$ at the point.

Many properties of these partially regular normal Cartan connections mirror those of regular normal connections. First of all, the lowest homogeneity component of its curvature is $\partial$-closed, meaning that it is still a section of the cohomology bundle
\be
\mc{P} \times_P H^2(\mf{g}^1,\mf{g}).
\ee

The first two BGG splitting operators are also defined for partially regular geometries. Furthermore, there is an invariant prolongation procedure for regular Cartan connections on any Tractor bundle $\mc{V}$. This paper demonstrates that this invariant prolongation procedure also works for the standard Tractor bundle $\mc{T}$, for most parabolic pairs.

The subsequent paper \cite{me2grad2} will then list and analyse all $|2|$-graded geometries, listing their cohomology spaces and their associated bundles, thus establishing when the above theorems apply. The case of rank six distributions on nine-dimensional manifolds will be analysed in full.

\subsection*{Acknowledgements}
It gives me great pleasure to acknowledge the financial support project P19500-N13 of the ``Fonds zur F\"orderung der wissenschaftlichen Forschung (FWF)'', as well as the help and comments of Andreas {\v{C}}ap and Jan Slov{\'a}k, and the support of the math department of the university of Vienna.

\section{The Cartan connection}
Let $G$ be a semisimple Lie group, and assume that there is a $|k|$-grading of its lie algebra
\be
\mf{g} = \mf{g}_{-k} \oplus \mf{g}_{-k+1} \oplus \ldots \oplus \mf{g}_{0 } \oplus \ldots \oplus \mf{g}_{k},
\ee
such that $[\mf{g}_i, \mf{g}_j] \subset \mf{g}_{i+j}$, $\mf{g}_0$ contains no simple summands of $\mf{g}$, and $\mf{g}_1$ generates all of $\mf{g}^1 = \sum_{j=1}^k \mf{g}_j$ via the Lie bracket. Define $\mf{g}^i$ as $\sum_{j= i }^k \mf{g}_{j}$. Then there is a corresponding filtration
\be
\mf{g} = \mf{g}^{-k} \supset \mf{g}^{-k+1} \supset \ldots \supset \mf{g}^{k-1} \supset \mf{g}^k = \mf{g}_k.
\ee
Let $P$ and $G_0$ be the subgroups of $G$ that preserve the filtration, respectively the grading, of $\mf{g}$. Then $P$ is a \emph{parabolic} subgroup of $G$.

Let $\mc{P} \to M$ be a principal $P$-bundle over a manifold $M$ of same dimensions as $\mf{g}^1$. A Cartan connection on $\mc{P}$ is a one-form $\omega$ on $\mc{P}$ taking values in $\mf{g}$, with the following properties:
\begin{enumerate}
\item $\omega$ is $P$-equivariant.
\item For any $A$ in $\mf{p}$, and $\xi_A$ the vector field on $\mc{P}$ corresponding to its action, $\omega (\xi_A) = A$.
\item For all $u \in \mc{P}$, $\omega_u : T\mc{P}_u \to \mf{g}$ is a linear isomorphism.
\end{enumerate}

The inclusion $P \subset G$ gives an extension of structure group $\mc{P} \subset \mc{G}$ and a unique principal connection $\omega'$ on $\mc{G}$ that is $G$-equivariant, has $\omega'(\xi_A) = A$ for all $A$ in $\mf{g}$, and pulls back to $\omega$ on $\mc{P}$. This descends to a linear connection on any vector bundle associated to $\mc{G}$, the so-called \emph{Tractor connection} $\onab$ (see \cite{TCPG}, or \cite{TBIPG}). Conversely, the Cartan connection $\omega$ can be derived from $\onab$ and the filtration of $\mc{A}$. The vector bundles associated to $\mc{G}$ are thus called Tractor bundles.

Now any representation $V$ of $G$ has a natural grading on it, coming from the grading of $\mf{g}$. We may thus talk of the homogeneity of any element of $V$. This grading is not preserved by the action of $P$; however, the corresponding filtration is. Consequently we have a filtration of the Tractor bundle $\mc{V} = \mc{P} \times_P V$. Thus we may talk about the \emph{minimal homogeneity} of any section of $\mc{V}$.

The explicit Tractor bundle we shall be using is the \emph{adjoint} Tractor bundle:
\be
\mc{A} = \mc{P} \times_P \mf{g} = \mc{G} \times_G \mf{g}.
\ee
Its natural filtration is:
\be
\mc{A} = \mc{A}^{-k} \supset \mc{A}^{-k+1} \supset \ldots \supset \mc{A}^{k},
\ee
where
\be
\mc{A}^j = \mc{P} \times_P \mf{g}^j.
\ee
Note that with the exception of $\mc{A}$ itself, none of these bundles are Tractor bundles. The minimal homogeneity of any section $\xi$ of $\mc{A}$ is now quite explicitly defined: it is the highest $j$ such that $\xi$ is also section of $\mc{A}^j$.

\subsection{The Levi and algebra\"ic brackets}
The Cartan connection gives an isomorphism $s_{\omega}$ between the tangent bundle $T$ and $\mc{A} / \mc{A}^0$. This isomorphism is called the \emph{soldering form}. The Killing form on $\mf{g}$ gives a metric on $\mc{A}$, and a consequent inclusion $T^* \subset \mc{A}$. Via $\tau$, $T$ inherits a filtration
\be
T = T^{-k} \supset T^{-k+1} \supset \ldots \supset T^{-1}.
\ee
We may define an associated graded bundle $gr(T)$, by $T_{j} = T_{j} / T_{j+1}$ and
\be
gr(T) = \sum_{j \geq -k}^{-1} T_{j}.
\ee
We may similarly define the graded $gr(\mc{A})$, and see that
\be
gr(\mc{A}) = gr(T) \oplus \A_0 \oplus gr(T^*).
\ee
Now $\mc{A}$ inherits an algebra\"ic bracket $\mc{K}$ from its Lie algebra structure. Its homogeneity zero component is invariant under the action of $P$. Consequently, $gr(\A)$ and $gr(T)$ will also inherit $\mc{K}$.

But $gr(T)$ also has another natural bracket, the \emph{Levi bracket} $\mc{L}$, coming from the Lie bracket. In detail, if $X$ and $Y$ are sections of $T_i$ and $T_j$, choose any lifts $\wt{X}$ and $\wt{Y}$ in $T^i$ and $T^j$, and define the section $\mc{L}(X,Y)$ of $T /T^{i+j+1}$ as
\be
\mc{L}(X,Y) = [\wt{X},\wt{Y}]/T^{i+j+1}.
\ee
The properties of the Lie bracket ensure that this is linear. Note that $\mc{L}$ is of non-positive homogeneity, and that the lowest homogeneity piece of $\mc{L}$ is well-define, independently of the choices of $\wt{X}$ and $\wt{Y}$. To define the higher homogeneities pieces of $\mc{L}$, we need an isomorphism $gr(T) \cong T$ to fix them. This will be accomplished in the next section, via the choice of a Weyl structure.

\subsection{Curvature of the Tractor connection}

The curvature $\kappa$ of $\onab$ is a section of $\wedge^2 T^* \otimes \A$. We may also see $\kappa$ as a $P$-equivariant map from $\mc{P} \to \wedge^2 \mf{g}_- \otimes \mf{g}$.

Dividing $\mc{P}$ by the action of $\exp \mf{g}^1$, there is a natural projection $\mc{P} \to \mc{G}_0$, where $\mc{G}_0$ is a principal $G_0$ bundle. We may choose a Weyl structures (see \cite{CapWeyl}) for $\onab$, see \cite{CapWeyl}. This is equivalently seen as a reduction of the structure group of $\A$ or $T$ from $\mc{P}$ to $\mc{G}_0$, an isomorphism $\A \cong gr(\A)$ (and consequent $T \cong gr(T)$), or a principal $\mc{G}_0$ connection $\nabla$ on $T$.

Given any Weyl structure, we have an explicit form for the Tractor connection. Given any vector field $X$ and any section $z$ of $\A$, it becomes:
\beqa \label{tractor:expression}
\onab_X z = \mc{K}(X,z) + \nabla_X z + \mc{K}(\rP(X),z),
\eeqa
where $\rP$ is a section of $T^* \otimes gr(T^*) \cong T^* \otimes T^*$. Furthermore, we now have a well-defined $\mc{L}: \wedge^2 gr(T) \to gr(T)$, at all non-positive homogeneities.

We are now in a position to explicitly calculate the curvature $\kappa$ in low homogeneities:
\begin{prop} \label{neg:curve}
In homogeneity $j <0$,
\be
\kappa_j = -\mc{L}_j,
\ee
the $j$-th homogeneity component of $\mc{L}$. In homogeneity zero
\be
\kappa_0 = \mc{K} - \mc{L}_0.
\ee
\end{prop}
\begin{proof}
Fix a Weyl structure $\nabla$. Note that in equation (\ref{tractor:expression}) for the Tractor connection, terms involving $\nabla$ or $\rP$ are of strictly positive homogeneity. For $X,Y$ sections of $T$ and $z$ a section of $\A$,
\be
\kappa(X,Y) z &=& \onab_X \onab_Y z - \onab_Y \onab_X z - \onab_{[X,Y]} z \\
&=& \mc{K}(X,\mc{K}(Y,z)) - \mc{K}(Y,\mc{K}(X,z)) - \mc{K}([X,Y],z) + \mathrm{ p.h.t.} \\
&=& \mc{K}( (\mc{K}(X,Y) - [X,Y] ) , z ) + \mathrm{ p.h.t.}
\ee
Here p.h.t.~designates positive homogeneity terms. The proposition then follows immediately from the expression for $\mc{L}$ and the fact that $\mc{K}$ is of homogeneity zero.
\end{proof}

Now, since $\nabla$ gives a reduction in structure group to $\mc{G}_0$, $\kappa$ may be seen as a $G_0$ equivariant map from $\mc{G}_0$ to $\wedge^2 \mf{g}^1 \otimes \mf{g}$. Designate by $C_j$ the space $\wedge^j \mf{g}^1 \otimes \mf{g}$. There are operators $\partial^* : C_j \to C_{j-1}$ and $\partial : C_j \to C_{j+1}$. They are adjoint operators for a certain definite signature metric of the $C_j$'s, see \cite{capslo}. The ones we will be using are $\partial^*$ on $C_2$ and $\partial$ on $C_1$. For $x$ and $y$ elements of $\mf{g}_- \cong (\mf{g}^1)^*$, the detailed expression is given as:
\beqa \label{partial:eq1}
(\partial \phi) (x,y) &=& \mc{K}(\phi(x),y) + \mc{K}(x,\phi(y)) - \phi(\mc{K}(x,y)) \\
(\partial^* \psi)(x) &=& \sum_l \mc{K}(z^l,\kappa(x,z_l)) - \frac{1}{2} \kappa(\mc{K}(z^l,x)_-,z_l), \label{partial:eq2}
\eeqa
for $(z_l)$ any basis of $\mf{g}_-$ and $(z^l)$ a dual basis of $\mf{g}^1$. Notice that these operators must preserve homogeneity (as $\mc{K}$ is of homogeneity zero), and are $G_0$-equivariant, so, given $\nabla$, extend to bundle operators on $\wedge^j gr(T) \otimes gr(\A)$. In fact, $\partial^*$ is $P$-equivariant, so it is well defined independently of the Weyl structure $\nabla$.

We will often need to split equation \ref{partial:eq2} into two components, $\partial^* =\alpha + \beta$, where
\beqa
\label{alpha} \alpha (\psi) (x) &=& \sum_l \mc{K}(z^l,\kappa(x,z_l)) \\
\label{beta} \beta (\psi) (x) &=& - \frac{1}{2} \sum_l \kappa(\mc{K}(z^l,x)_-,z_l).
\eeqa

\begin{defi}[Normality]
A Cartan connection $\omega$ is normal if and only if the curvature $\kappa$ of the corresponding Tractor connection $\onab$ is $\partial^*$-closed, i.e.
\be
\partial^* \kappa = 0.
\ee
\end{defi}

Since $(\partial)^2 = (\partial^*)^2 = 0$, there are cohomology groups $H^j(\mf{g}^1, \mf{g})$ for $\partial^*$, defined as $ker \ \partial^*_j / \ im  \ \partial^*_{j+1}$. Since $\partial$ and $\partial^*$ are adjoint for a definite signature metric (\cite{capslo}), this space is isomorphic with $ker \ \partial_j / \ im  \ \partial_{j-1}$ and can be identified as the kernel of the Beltrami-Laplacian
\be
\square = \partial \partial^* + \partial^* \partial.
\ee
There is further Hodge decomposition of $C_j$ as
\beqa \nonumber
C_j &=& im \ \partial^* \oplus im \ \partial \oplus H^j(\mf{g}^1,\mf{g}) \\
\label{hodge} &=& im \ \partial^* \oplus im \ \partial \oplus ker \ \square.
\eeqa

\subsection{$|2|$-graded geometries}

If $\mf{g}$ is $|2|$-graded, then minimum homogeneity of $\wedge^2 T^* \otimes \A$ is zero; hence $\mc{L}$ is of homogeneity zero, and hence is, like $\mc{K}$, a section of $\wedge^2 T^*_1 \otimes T_{-2}$.

\begin{defi}[$|2|$-graded pre-geometry]
A $|2|$-graded pre-geometry is a manifold $M$ with tangent bundle $T$ and a distribution $H \subset T$ and a $|2|$-graded parabolic pair $(G,P)$ such that
\be
\textrm{rank} \textrm{ } (H) &=& \textrm{dim} \textrm{ } \mf{g}_{-1} \\
\textrm{rank} \textrm{ } (H/T) &=& \textrm{dim} \textrm{ } \mf{g}_{-2}.
\ee
\end{defi}

The distribution $H$ is the candidate for being $T_{-1}$, while $H/T$ is the candidate for $T_{-2}$.

We are now ready to define a partially regular geometry:
\begin{defi}[Partially regular]
A partially regular $|2|$-graded geometry is a $|2|$-graded pre-geometry along with a principal $P$-bundle $\mc{P} \to M$ that is the structure bundle for $gr(T) = H \oplus T/H$.

This structure bundle is enough to define the algebra\"ic bracket $\mc{K}$ on $gr(T)$ and the $\partial^*$-operator in homogeneity zero. Then the geometry is partially regular if
\be
\partial^* \mc{K} - \mc{L} = 0.
\ee
\end{defi}

Similarly, the definition of a partially regularisable pre-geometry is the obvious one:
\begin{defi}[Partially regularisable]
A partially regularisable $|2|$-graded geometry is a $|2|$-graded pre-geometry such that there exists a principal $P$-bundle $\mc{P} \to M$ for $gr(T)$ making $M$ into a partially regular geometry.
\end{defi}

\section{Normalising homogeneity zero}

If we have a manifold $M$ with a distribution $H \subset TM$, then let $\mc{S}$ be the graded principle bundle for the graded bundle $gr(TM) = H \oplus TM/H$. Thus $\mc{S}$ is a $S = GL(m_1) \times GL(m_2)$-bundle, with $m_1$ the rank of $H$ and $m_2$ the rank of $TM/H$.

Note that the Levi bracket $\mc{L}$ is a section of
\be
\mc{S} \times_S \wedge^2 \mbb{R}^{{m_1}} \otimes \mbb{R}^{m_2}.
\ee

We may make the same definitions for $M$, $H$, $\mc{S}$ and $S$ in the complex case. Then if we furthermore require that $H$ be a holomorphic distribution, we get the similar result that $\mc{L}$ must be a section of
\be
\mc{S} \times_S \wedge^2 \mbb{C}^{{m_1}*} \otimes_{\mbb{C}} \mbb{C}^{m_2}.
\ee
We will use the same notation in the real and holomorphic category, distinguishing between the two by context alone.

\begin{defi}[Open, maximum measure]
Given a manifold $M$, we will call a subset $M' \subset M$ open, maximum measure (OMM) if $M'$ is open and the complement $M - M'$ is of strictly lower dimension than $M$.
\end{defi}

The main results then are:
\begin{theo} \label{main:zero:theo}
Let $M$ be a complex manifold with tangent space $T$, let $H \subset T$ be a holomorphic distribution, and let $(G,P)$ be a complex parabolic $|2|$-graded pair such that
\be
\textrm{rank}_{\mbb{C}} \textrm{ } (H) &=& \textrm{dim}_{\mbb{C}} \textrm{ } \mf{g}_{-1} \\
\textrm{rank}_{\mbb{C}} \textrm{ } (H/T) &=& \textrm{dim}_{\mbb{C}} \textrm{ } \mf{g}_{-2}.
\ee
This makes $(M,H,G,P)$ into a holomorphic $|2|$-graded pre-geometry.

Let $W = \wedge^2 \mf{g}_{-1}^* \otimes \mf{g}_{-2}$, where the tensor and wedge products are complex. Then there exists a set $U \subset W$, Zariski-open in $W$ and $S$-invariant. This $U$ has the property that if there is a point $x \in M$ such that $\mc{L}_x$ is in
\be
\mc{S}_x \times_S U,
\ee
then there exists an OMM submanifold $M' \subset M$ such that $H|_{M'}$ is partially regularisable. Moreover, if $H^1(\mf{g}^1,\mf{g})$ vanishes in homogeneity zero, this partial regularisation is unique up to $G_0$ action and a finite choice.
\end{theo}
Zariski-open means OMM, of course. Colloquially, the previous Theorem states that `almost all holomorphic distribution of the correct rank and co-rank are partially regularisable, almost everywhere'. The `finite choice' aspect of the Theorem is unsatisfactory, however:
\begin{theo} \label{main:zero:theo:real}
In the conditions described above, let $u_P$ be a uniqueness procedure on $W$ (see section \ref{uniq:proc}). Then there exists a subset $U' \subset U \subset W$, OMM in $W$ and dependent on the choice of $u_P$. This $U'$ has the property that if there is a point $x \in M$ such that $\mc{L}_x$ is in
\be
\mc{S}_x \times_S U',
\ee
then there exists an OMM $M' \subset M$ such that $H|_{M'}$ is partially regularisable. Moreover, if $H^1(\mf{g}^1,\mf{g})$ vanishes in homogeneity zero, this partial regularisation is unique up to $G_0$.
\end{theo}
The uniqueness procedure will be described in more detail in section \ref{uniq:proc}. Simply put, the problem is that points in $W$ on different $G_0$ orbits, corresponding to different partial regularisations, might be on the same $S$-orbit. The uniqueness procedure then chooses a single $G_0$ orbit representative.

The above is true in the holomorphic category, but what is true in the real category? Unfortunately, the proof relies on algebra\"ic geometry results that do not carry through to the real case. Though I strongly believe the more general result is true in the real case, the best that can be proved is:
\begin{theo}
Let $M$ be a real manifold with tangent bundle $T$, let $H \subset T$ be a distribution, and let $(G,P)$ be a parabolic $|2|$-graded pair such that
\be
\textrm{rank } (H) &=& \textrm{dim } \mf{g}_{-1} \\
\textrm{rank } (H/T) &=& \textrm{dim } \mf{g}_{-2}.
\ee
This makes $(M,H,G,P)$ into a $|2|$-graded pre-geometry.

Let $W = \wedge^2 \mf{g}_{-1}^* \otimes \mf{g}_{-2}$. Then there exists an open set $U \subset W$, $S$-invariant. This $U$ has the property that if there is a point $x \in M$ such that $\mc{L}_x$ is in
\be
\mc{S}_x \times_S U,
\ee
then there exists an open submanifold $M' \subset M$ containing $x$ such that $H|_{M'}$ is partially regularisable. Moreover, if $H^1(\mf{g}^1,\mf{g})$ vanishes in homogeneity zero, this partial regularisation is unique up to $G_0$ action and a finite choice.

If we have uniqueness procedure $u_P$, then the same result holds as above, giving $U' \subset U \subset W$ OMM in $U$ and a unique choice of partial regularisation.
\end{theo}
The result is somewhat stronger (and does not need a uniqueness procedure) if we start with a partially regularisable distribution and then perturb it; though we will need to wait till the end of this section in order to phrase it precisely.

The rest of this section will be devoted to proving the above Theorems.

Now the Levi bracket $\mc{L}$ is equivalent with an $S$-equivariant function $f_{\mc{L}}: \mc{S} \to W$. If $\mc{K}$ is the standard bracket in $W$, define $V$ as the $S$-orbit of $\mc{K}$ in $W$. Then a choice of algebra\"ic bracket $\mc{K}$ is a choice of $S$-equivariant function $f_{\mc{L}}: \mc{S} \to V$. At each point $x \in M$, $\partial^*_{\mc{K}} \mc{K} - \mc{L} = 0$ if and only if
\beqa \label{partial:func}
\partial^*_{f_{\mc{K}}(u)} f_{\mc{K}}(u) - f_{\mc{L}}(u) = 0,
\eeqa
where $u$ is some point in $\mc{S}_x$.

Set $k = f_{\mc{K}}(u)$ and $l = f_{\mc{L}}(u)$. Then to define $\partial^*_{k}$, we need the element $k^* \in W^*$. This $k^*$ may be defined by noting that the stabiliser of $k$ in $W$ is conjugate to $G_0$, and preserves a line in $W^*$. Then $k^*$ is the element of that line whose complete contraction with $k$ is unity.

We will switch to abstract index notation, and denote $k$ by $k_{ab}^A$, with the small latin indices denoting $\mf{g}_{-1}^*$, and the capital latin indexes denoting $\mf{g}_{-2}$. Then equation (\ref{partial:func}) can be rewritten as:
\beqa \label{k:one}
(k^*)_{A}^{ab}(k _{ac}^A - l_{AB}^A) &=& 0 \\
\label{k:two} -\frac{1}{2}(k^*)_A^{ab}(k_{ab}^B - l_{ab}^B) &=& 0.
\eeqa
These two equations are the same equations for $\alpha$ and $\beta$ as (\ref{alpha}) and (\ref{beta}). In homogeneity zero, $\alpha$ maps to $End(\mf{g}_{-1})$ and $\beta$ maps to $End(\mf{g}_{-2})$; consequently, if a quantity is $\partial^*$-closed, then it must be $\alpha$-closed and $\beta$-closed.

Note that the above equations are invariant under the action of $S$, so do not depend on the choice of $u \in \mc{S}_x$.

The process to partially regularize $H$ at $x$ would then be to construct a $k$ at $u$ that solves equations (\ref{k:one}) and (\ref{k:two}), and then extend to all of $\mc{S}_x$ by $S$ action, giving $\mc{K}_s$. One could imagine starting with a random $k \in V$, and then acting on it by $S$, until one finds a candidate solving the equations.

In fact, it is conceptually simpler to fix $k$ and then act on $u$ via $S$. Since $f_{\mc{L}}$ is $S$-equivariant, this is equivalent with picking $l$ and acting on it with $S$. So the values of $l$ (more precisely, the $S$-orbit of $l$) is what determines whether $H$ is partially regularisable at $x$. Since $k$ is the standard bracket, we will drop the subscript from $\partial^*_k$ and just write it $\partial^*$.

This problem can be treated through algebra\"ic geometry. Now $S = S_1 \times S_2$, where $S_i = Aut(\mf{g}_{-i})$. We define $\overline{S}$ to be the affine space $End(\mf{g}_{-1}) \oplus End(\mf{g}_{-2})$. Then we define the affine variety $\mc{W}$ via: 
\be
\mc{W} = \{(w,s) \in W \times \overline{S} | \partial^* (k - s \cdot w) = 0 \},
\ee
where
\be
(s \cdot w)(x,y) = s_2(w(s_1(x),s_1(y)))
\ee
for $x, y \in \mf{g}_{-1}$. This is just an extension of the group action of $S$ to the affine space $\overline{S}$. Moreover, this is a polynomial action, while $\partial^*$ is linear, hence the whole equation defines an algebra\"ic variety -- non-empty, since $(k,Id)$ is certainly an element of $\mc{W}$. Then define $p: \mc{W} \to W$ as the projection onto the first component; then we can partially regularize $H$ at $x$ if and only if $l \in p(\mc{W})$.

\begin{prop} \label{working:prop}
In the holomorphic category, the set $p(\mc{W})$ contains a subset $U$ that is Zariski open in $W$. If $H^1(\mf{g}^1,\mf{g})$ is zero in homogeneity zero, then the fibers of $p$ at any point of $U$ are isomorphic to a finite collection of right cosets of $G_0$.
\end{prop}
\begin{proof}
The main body of the proof will come from the following proposition, a simplified version of Proposition 8.8.1 and Corollary 8.8.2 from the book \cite{algeosavior}:
\begin{prop} \label{al:geo}
Let $f$ be a polynomial map between complex quasi-affine algebra\"ic varieties. Then the image of $f$ contains a Zariski open subset of its Zariski closure, on which the fibers of $f$ are of minimal dimension (i.e.~the dimension of the fibers of $f$ plus the dimension of the image, is equal to the dimension of the domain).
\end{prop}
Now define
\be
\Theta : W \times \overline{S} &\to& \textrm{ im } \partial^* \\
\Theta(w,s) &=& \partial^* (k - s \cdot w).
\ee
The variety $\mc{W}$ is therefore the zero set of $\Theta$. We want to take the derivative of $\Theta$ in the $\overline{S}$ directions, around the point $(w,Id)$; then it is easy to see from the formula of the action of $\overline{S}$ on $W$, that:
\beqa \nonumber
D\Theta(w,Id)(t_1,t_2) &=& \partial^*( - t_2(w(-,-)) -w(t_1(-),-) - w(-, t_1(-))) \\
\label{deriv:general} &=& -\partial^* (\partial_w (t_1+t_2) ),
\eeqa
where $t_i \in \mf{s}_i$. Notice that $\mf{s}_i = \overline{S}_i$ as spaces, but we use a different notation to keep track of the different actions: an algebra action for $\mf{s}_i$, and the extension of a group action for $\overline{S}_i$. Now if we substitute $k = w$, the expression becomes
\beqa \label{deriv:form}
D \Theta(k,Id) = -\partial^* \partial.
\eeqa
The rest of the proof is by dimension counting. We know that the space $\mf{s}_1 \oplus \mf{s}_2$ decomposes, via equation (\ref{hodge}), into the sum
\be
\textrm{im } \partial^* \oplus H^1(\mf{g}^1,\mf{g})_0 \oplus \textrm{im } \partial.
\ee
By equation (\ref{deriv:form}) and the properties of $\partial$ and $\partial^*$, the fibre of $p$ at $k$ is of same dimension as ker $\partial = H^1(\mf{g}^1,\mf{g})_0 \ \oplus$ im $\partial$. However, the restriction $\Theta(w,s) = 0$ implies that the codimension of $\mc{W}$ in $W \times \mc{S}$ is the same as the dimension of ker $\partial$. Thus the dimension of the image of $p$ around $k$ is equal to the dimension of $W \times \overline{S}$, minus the dimension of $\overline{S}$ -- in other words, it is equal to the dimension of $W$.

Thus $p(\mc{W})$ contains a Euclidean open set around $k$. The Zariski closure of a Euclidean open set must be the whole of $W$, so Proposition \ref{al:geo} demonstrates the existence of $U$, and furthermore that the fibers of $p$ on $U$ are all of constant dimension, equal to the dimension of ker $\partial$. There must be finitely many such fibers, as each one is an algebra\"ic variety with a bounded degree, and thus with a bounded number of components.

An extra subtlety must be noted: it is possible that some elements of $U$ come only from elements $(w,s)$ where $s \in \overline{S}$ is not invertible. In general, this will not be the case (as $\partial^* k$ in invertible), and we may avoid the problem entirely if we restrict to the \emph{quasi-affine} variety:
\be
\mc{W} \cap \{(w,s)| \det s \neq 0\}.
\ee

If $H^1(\mf{g}^1,\mf{g})_0 = 0$, then the fibers are each of same dimension as im $\partial$ -- hence of same dimension as $G_0$. We then merely need to note that for $g \in G_0$, if $\Theta(w,s) = 0$, then
\be
\Theta(w,gs) &=& \partial^* (k - (g \cdot s \cdot w)) \\
&=& \partial^* ((g \cdot k) - (g \cdot s \cdot w)) \\
&=& g \cdot \partial^* (k - s \cdot w) \\
&=& g \cdot \Theta(w,s) = g \cdot 0 = 0,
\ee
since $G_0$ fixes $k$ and commutes with $\partial^*$. Thus $S$ components of the fibers of $p$ on $U$ must consist precisely of a finite collection of left cosets of $G_0$ in $\overline{S}$.

\end{proof}

Now we may prove the general theory in the holomorphic category. The $U \subset W$ has been defined above, and is evidently $S$-invariant. If $f_{\mc{L}}(u) \in U$, then $H$ is partially regularisable at $x$. Since $U$ is open and $\mc{L}$ is continuous, the $H$ is partially regularisable in a neighbourhood of $x$. Finally, since $W-U$ is an algebra\"ic subvariety, and $H$ is holomorphic, the subset $N$ of $M$ where $f_{\mc{L}} \notin U$ is an analytic subvariety, hence $M' = M - N$ is open, dense, and full measure in $M$.

In the real category, we do not have such general results as Proposition \ref{al:geo}; so all we can claim is that $p(\mc{W})$ has non-empty interior, and hence that $U$ can be chosen to be open. Consequently $M'$ is also open.

The result in the real case is disappointing, and even in the holomorphic case, the `finite choice' may be very high. This can be fixed by setting a uniqueness procedure see section \ref{uniq:proc}. But there is an elegant method for guaranteeing uniqueness when we are perturbing the right sort of distribution. First, we shall define what a `suitable' partial regularisation is:
\begin{defi}
Let $M$ be a manifold and $H \subset TM$ a distribution of correct rank and co-rank, as above. Assume that for $x \in M$, there exists a partial regularisation $\mc{K}_x$ -- equivalently, a reduction in the structure group of $gr(T)_x$ to $(\mc{G}_0)_x$. Then this partial regularisation is \emph{suitable} if for a given $u \in (\mc{G}_0)_x$, we have $l = f_{\mc{L}}(u)$ such that the $S$-orbit of $l$ is transverse to the kernel of $\theta(w) = \partial^* (k - w)$.
\end{defi}
Since the above relation is $G_0$ invariant, it does not depend on the choice of $u \in (\mc{G}_0)_x$. In practice, finding out if a given regularisation is suitable is very tricky; however, there is a generic sufficient condition that is much easier to calculate. Namely that $H$ is suitably partially regularised at $x$ if
\be
\partial^* \partial_{l}
\ee
is a maximum rank map from $\mf{s}$ to im $\partial^*$ (recall that this is the derivative of $\Theta$ at $(l,Id)$ in the $S$ directions). The above is evident by implicit function theorem; and equation (\ref{deriv:form}) implies the immediate corollary that:
\begin{cor}
A regular geometry is always a suitable partial regularisation, at every point.
\end{cor}

The main result comes from the following lemma:
\begin{lemm}
Let $H_t$, $t \in (a,b) \subset \mbb{R}$ be a continuous family of constant rank distributions. Assume there is a point $x \in M$, a $\tau \in (a,b)$ and a given suitable partial regularisation of $H_{\tau}$ at $x$, and that $H^1(\mf{g}^1,\mf{g})_0 =0$. Then there exists a neighbourhood $N$ of $x$ in $M$ and an open interval $I\subset \mbb{R}$ containing $\tau$ such that there are unique suitable partial regularisations of $H_t$, on the set $N$, for $t \in(\tau - \epsilon, \tau + \epsilon)$, continuous in $t$ and on $N$.
\end{lemm}
\begin{lproof}
Choose any neighbourhood $N$ on which $\mc{S}$ is trivial, fix a given $u \in (\mc{G}_0)_x$, and extend $u$ to a local section $\mu$ of $\mc{S}|_N$. Let $I \subset \mbb{R}$ be an open interval containing $\tau$. Then we may see $f_{\mc{L}}$ as a function on $\mc{S}|_N \times I$; then $f_{\mc{L}} \circ \mu$ is a map from $N \times I$ to $W$.

Now our partial regularisation of $H_{\tau}$ at $x$ is suitable, so the $S$ orbit of $l = f_{\mc{L}}(u) = f_{\mc{L}} \circ \mu(x)$ is transverse to the kernel of $\theta$. This is a open condition, so remains true on a small neighbourhood of $l$. Let us then look at the subset $Id \ \oplus$ im $\partial$. Since this is transverse to $G_0$ at $Id$, and $G_0$ preserves the kernel of $\theta$, we know that locally around $l$, $W$ decomposes as
\be
A \times \Lambda,
\ee
where $\theta(A) = 0$ and $\Lambda \subset (Id \ \oplus$ im $\partial)$. This decomposition is according to the orbits of points in $A$ under the action of $\Lambda$.

Then, restricting to open subsets of $N$ and $I$ as needed, the above implies there is a unique choice of $g: N \times I \to S$ such that $g \cdot \mu$ defines a partial regularisation of $H_{t}$ on $N$ for $t \in I$. This partial regularisation is seen to be suitable, by the above decomposition of $W$.

Changing our choice of $\mu$ just changes this partial regularisation by the action of $G_0$, hence our definition is unique.
\end{lproof}
An immediate consequence of the above result is that:
\begin{prop}
Let $M$ be a manifold and let $H_t$, $t\in (a,b) \in\mbb{R}$ be a continuous family of distributions on $M$, such that $H_t$ is independent of $t$ outside of a compact subset $M' \subset M$. Assume that $H_{\tau}$ is suitably partially regularised and $H^1(\mf{g}^1,\mf{g})_0 =0$. Then there exists an open interval $I \subset \mbb{R}$ containing $\tau$ such that there is a unique continuous suitable partial regularisation of $H_t$ for $t \in I$.
\end{prop}
\begin{proof}
By the previous lemma, for each point $x$ in $M'$, there exists neighbourhood $N_x$ and $I_x$ such that we can extend the suitable partial regularisation uniquely on $N_x \times I_x$. The $N_x$ form a cover of $M'$; since $M'$ is compact, there exists a finite set $F \subset M$ such that $\{N_x| x \in F\}$ is a cover of $M$. Then set $I = \cap_{x \in F} I_x$, and we have the proposition.
\end{proof}
\begin{rem}
The above means that for small deviations from regular geometries on compact sets, we have unique partial regularisations.
\end{rem}
In general, we will only get uniqueness through some type of uniqueness procedure; this will be the talk of the next section.

\subsubsection{Uniqueness procedure} \label{uniq:proc}
By the previous results, we know that there exists a set $U$, Zariski open in $W$ and closed under the orbit of $S$, such that for each $l\in U$, the orbit space $S\cdot l$ meets the kernel of $\theta$ in a finite collection of $G_0$-orbits. The number of these $G_0$ orbits is bounded for varying $l$, by considerations of degree. Denote by $K$ the kernel of $\theta$. By picking the standard bracket $k$ as the origin in $K$, we may see $K$ as a vector space, not only as an affine space.

It would be ideal to have a method for selecting a single $G_0$ orbits from among these finite collections. Such a method may not be valid on the whole of $U$; if it valid on an OMM $U' \subset U$, we call it a uniqueness procedure $u_P$.

For a given $l \in W$, call $d_l$ the dimension of the stabilizer subgroup of $l$ in $S$ (since the grading element is in $S$ and fixes all of $W$, $d_l \geq 1$). If $d_l > d_{l'}$ for $l,l' \in W$, then almost all elements $m$ of the affine space spanned by $l$ and $l'$ will have $d_{m} \leq d_l$. Thus there exists a $d$ such that $d \leq d_l$ for all $l \in W$ and $d = d_l$ an OMM subset of $W$. This set must intersect $U$ in another OMM subset.

Now pick an $l \in U$ such that $d_l = d$ and $l \in K$ (since the $S$ orbit of any element in $U$ intersects $K$, this can always be done). Then pick a hermitian metric $h$ on $K$.

\begin{prop}
The element $l$ along with $h$ defines a uniqueness procedure.
\end{prop}
\begin{proof}
Since $l \in K$, the stabilizer subgroup of $l$ must be contained in $G_0$. Since $W$ is a vector space, we may identify a subundle $K_l \subset TW_l$, isomorphic with $K$. Let $D_l : \mf{s} \to TW_l$ be the derivative of the $S$ action on $l$, and define the subspace $\Lambda \subset K$ as $\Lambda = ((D_l (\mf{g}_0))^{\perp}) \cap K_l$. Since $W$ is a vector space, we may see $\Lambda$ as an affine subspace through $l$, contained in $K$.

By construction, $\Lambda$ is transverse to the $S$-orbit of $l$ (since the action of $S/G_0$ maps $l$ off of $K$). We may define the algebra\"ic variety
\be
\mc{W}^{\lambda} = \{(w,s) \in W \times \overline{S} | s \cdot w \in \Lambda \}.
\ee
Similarly to the proof of Proposition \ref{working:prop}, we can see that $p(\mc{W})$ has non-empty interior. Hence by Proposition \ref{al:geo} there is a OMM subset $U' \subset U$ such that $U'$ is $S$-closed, and for every $l' \in U'$, the orbit $S \cdot l'$ meets $\Lambda$ transversally. By dimension count, this means that $S \cdot l '$ meets $\Lambda$ only in isolated points; considerations of degree imply that the number of isolated points is finite. Moreover, for any two points in $U'$, we can consider the orbit of the affine space generated by those two points, and its intersection with $\Lambda$; this demonstrates that for an OMM subset $U'' \subset U'$, the cardinality of this intersection is a constant across $U''$.

Then, whenever it is possible, we pick a representative of $S \cdot l'$ by choosing the point $s\cdot l'$ in $(S \cdot l') \cap \Lambda$ such that $||s\cdot l' - l||$ is minimized -- using $h$ to calculate this norm.

This procedure is evidently well defined and continuous around $l$. Then note that the norms of the different intersection points of $(S \cdot l') \cap \Lambda$ are locally real analytic functions on $U''$; hence the above procedure is well-defined, and continuous, for a OMM subset $U''' \subset U''$.

The uniqueness procedure is then finalised by selecting the $G_0$-orbit of the chosen representative.
\end{proof}

\begin{rem}
Ideally, we would want to pick an $l$ close to the standard bracket $k$; that way, the uniqueness procedure extends continuously to $k$ (the chosen representatives of points near $k$ need not converge to $k$; but there will exist points in the $G_0$ orbits of the chosen representatives that converge to $k$).
\end{rem}
\begin{rem}
Note that the final choice of $U'''$ involves a real function (the metric), and hence $U'''$ is not a complex pseudo-affine variety, but a real one -- in particular, it may be disconnected.
\end{rem}
\begin{rem}
The above procedure generalises to the real category, with the proviso that $U'''$ is open, but need no longer be OMM.
\end{rem}
The dependence on $l$ is somewhat problematic, so I conjecture that:
\begin{con}
The choice of $h$ alone is enough to define a uniqueness procedure. The $\Lambda$ is then defined by choosing a generic sequence of $l_j$ with $d_{l_j} = d$ and $l_j \to k$. Then if $\Lambda(l_j)$ tend to a well defined limit $\Lambda$ of same dimension, and if this $\Lambda$ is independent of the choice of generic sequence $l_j$, then we may use the same procedure as above, substituting $k$ for $l$.
\end{con}

\subsection{Non-regular path geometries}
Paper \cite{meCR} demonstrates that there can be an analogue of `partial regularisation', even for structures where $H^1(\mf{g}^1,\mf{g})$ does not vanish in homogeneity zero -- in this case, co-dimension one CR structures. The redundancy provided by $H^1(\mf{g}^1,\mf{g})_0$ is removed by fixing a choice of extra structure. For the co-dimension one CR structure, this extra structure is a complex structure.

There is an analogue construction for path geometries. The extra piece of information is provided by the
\begin{defi}[Non-regular path geometry]
A non-regular path geometry is given by a manifold with a distribution $H$ of rank $n+1$ and co-rank $n$, and a line-bundle $L \subset H$.
\end{defi}
With the extra piece of information given by $L$, non-regular path geometries can be treated identically to other non-regular geometries:
\begin{theo}
All the theorems of this section apply to non-regular path geometries, with the group $S$ restricted to the subgroup $A \subset S$ that preserves the distinguished line.
\end{theo}
\begin{proof}
For complex path geometries, the algebra components $\mf{g}_{-1}$ and $\mf{g}_{-2}$ are
\be
\mf{g}_{-1} &=& \mbb{C}^{n*} \oplus \mbb{C} \\
\mf{g}_{-2} &=& \mbb{C}^{n*}.
\ee
The space $H^1(\mf{g}^1, \mf{g})_0$ can be calculated, by Kostant's proof of the Bott-Borel-Weyl theorem (see \cite{Kostant}); it is $\mbb{C}^n$, and lies inside $\mf{gl}(\mf{g}_{-1})$ as the strictly diagonal algebra mapping $\mbb{C}$ to $\mbb{C}^{n*}$. Hence note that if $\mf{a}$ is the Lie algebra of $A$,
\be
\mf{s} = \mf{a} \oplus H^1(\mf{g}^1, \mf{g})_0.
\ee
Thus if we restrict to $A$ (an algebra\"ic subset of $S$), all the results of the previous section apply, including the uniqueness results, since $\mf{a} \cap H^1(\mf{g}^1, \mf{g})_0 = 0$.
\end{proof}

\begin{rem}
This result can be extended to more general Grassmannian geometries, but the picture becomes more complicated. The group $A$ must be replaced by a product of groups, and the definition of a non-regular Grassmannian geometry is no longer given by a single distinguished subbundle of $H$, but by a class of such subbundles.
\end{rem}

\section{Normalising higher homogeneities}

\begin{theo} \label{theo:high}
Given a partially regular $|2|$-graded geometry $(M,H,\mc{P})$, there exists a partially regular \emph{normal} Tractor connection $\onab$ on $\mc{P}$. This Tractor connection is unique up to automorphism if $H^1(\mf{g}^1,\mf{g})$ is zero in strictly positive homogeneities.
\end{theo}

Note that the cohomology condition is satisfied for all parabolic pairs $(G,P)$ that do not have a summand of projective type
\be
\mf{sl}(n+1, \mbb{C}) \ / \ \mf{gl}(n,\mbb{C}) \rtimes \mbb{C}^n
\ee
(or any real forms of the above) or of contact projective type 
\be
\mf{sp}(2n+2,\mbb{C}) \ / \ (\mbb{C} \oplus \mf{sp}(2n,\mbb{C})) \rtimes \mbb{C}^{2n} \rtimes \mbb{C}
\ee
(or any real forms of the above).

This theorem will be an immediate consequence of the following proposition:
\begin{prop} \label{curvature:adjustment}
Let $(\mc{P},\onab)$ be any partially regular $|2|$-graded Tractor connection, and let $\theta \in T^* \otimes \A$ be of homogeneity $l \geq 1$. Let $\kappa$ be the curvature of $\onab$ and $\kappa_{\theta}$ the curvature of $\onab+ \theta$. Then
\begin{itemize}
\item $\kappa_{\theta} - \kappa$ is of homogeneity $ \geq l$,
\item $(\kappa_{\theta} - \kappa)_l = \Psi(\theta_l)$, where $\Psi$ is a linear map from $(T^* \otimes \A)_l$ to $(\wedge^2 T^* \otimes \A)_l$,
\item using the soldering form $s_{\theta}$ of $\onab + \theta$ to identify $T$ and $\A / \A^0$, $\partial^* \circ \Psi$ is a linear map from $(T^* \otimes \A)_l$ to itself, which is invertible on the image of $\partial^*$.
\end{itemize}
\end{prop}

\begin{proof}
For the rest of this proof, let $X$ and $Y$ be sections of $\A / \A^0$, and $X'$ and $Y'$ be sections of $gr(\A / \A^0)$.

Let $s$ be the soldering form for $\onab$, and $s_{\theta}$ the soldering form for $\onab + \theta$. Then it is easy to see that $s_{\theta} = s + \hat{\theta}$, where $\hat{\theta}$ is the projection of $\theta$ onto $T^* \otimes (\A / \A^0)$. Since $\theta$ is of homogeneity $l > 0$, $s_{\theta}$ remains invertible (since its homogeneity zero piece is invertible), and in homogeneity $l$,
\be
s_{\theta}^{-1} = s^{-1} - s^{-1} \circ \hat{\theta} \circ s^{-1}.
\ee
Then define $\widetilde{\theta} = \hat{\theta} \circ s^{-1}$. By a slight abuse of notation, we will also designate $\theta \circ s^{-1}$ by $\widetilde{\theta}$, and, similarly for the curvature $\kappa$:
\be
\widetilde{\kappa}(X,Y) = \kappa(s^{-1}(X),s^{-1}(Y)).
\ee

Then the curvature of $\onab + \theta$, seen, via $s_{\theta}$, as a section of $\wedge^2 (\A^{1}) \otimes \A$, is
\be
\kappa_{\theta}(s_{\theta}^{-1}(X),s_{\theta}^{-1}(Y)) &=& \kappa_{\theta}(s_{\theta}^{-1}(X),s_{\theta}^{-1}(Y)) + \onab_{s_{\theta}^{-1}(X)} \theta(s_{\theta}^{-1}(Y)) \\
&& - \onab_{s_{\theta}^{-1}(Y)} \theta(s_{\theta}^{-1}(X)) - \theta([s_{\theta}^{-1}(X),s_{\theta}^{-1}(Y)]) \\
&=& \widetilde{\kappa}_l(X,Y) - \widetilde{\kappa}_0(\widetilde{\theta}(X),Y) - \widetilde{\kappa}_0(X, \widetilde{\theta}(Y)) \\
&& + \mc{K}(X,\widetilde{\theta}(Y)) - \mc{K}(Y,\widetilde{\theta}(X)) \\
&& - \widetilde{\theta}(\mc{K}(X,Y)) + \widetilde{\theta} (\widetilde{\kappa}_0(X,Y)),
\ee
plus terms of homogeneity $\geq l$. This demonstrates that $\kappa - \kappa_{\theta}$ is of homogeneity $l$. We may rewrite this more succinctly as saying, in homogeneity $l$, that
\be
\big(\kappa_{\theta}(s^{-1}_{\theta}(X'), s^{-1}_{\theta}(Y')) - \widetilde{\kappa}_l(X',Y') \big)_l = \partial(\widetilde{\theta}_l)(X',Y') + \widetilde{\theta}_l \bullet \widetilde{\kappa}_0(X',Y'),
\ee
where $\widetilde{\theta}_l \bullet \widetilde{\kappa}_0(X',Y') = \widetilde{\theta}_l (\widetilde{\kappa}_0(X',Y')) - \widetilde{\kappa}_0(\widetilde{\theta}_l(X'),Y')) - \widetilde{\kappa}_0(X', \widetilde{\theta}_l(Y')) $. Hence
\be
\Psi({\theta})(s^{-1}_{\theta} , s^{-1}_{\theta} )_l  = \partial(\widetilde{\theta}_l) + \widetilde{\theta}_l \bullet \widetilde{\kappa}_0.
\ee
We will define $\widetilde{\Psi}: (\A^1) \otimes \A \to \wedge^2 (\A^1) \otimes \A$ by saying the above expression is equal to $\widetilde{\Psi}(\wt{\theta}_l)$. It now remains to prove the rest of the proposition, namely that:

\begin{lemm}
The map $\partial^* \circ \wt{\Psi}$ is invertible on the image of $\partial^*$.
\end{lemm}
\begin{lproof}
First note that
\be
\partial^* \circ \wt{\Psi} (\wt{\theta}_l) = \square \wt{\theta}_l + \partial^* (\widetilde{\theta}_l \bullet \widetilde{\kappa}_0).
\ee
To demonstrate that this is invertible on the image of $\partial^*$, we shall show that $\partial^* \circ \wt{\Psi}$ is `upper-triangular' in some sense, and that $\partial^* (\widetilde{\theta}_l \bullet \widetilde{\kappa}_0)$ is the strictly upper-triangular piece. Then invertability will flow from the invertability of the Beltrami-Laplacian $\square$, see \cite{capslo}.

First note that $\square$ respects the decomposition 
\be
(\A / \A^0 \otimes \A)_l = \big(\A_1 \otimes \A_{l-1}\big) \oplus \big( \A_2 \otimes \A_{l-2} \big).
\ee

Now assume that $\wt{\theta}_l$ is a section of $\A_1 \otimes \A_{l-1}$. Then $\wt{\theta}_l(X')$ is a section of $\A_{l-1}$, and, since $l >0$,
\be
\big( \wt{\theta}_l(X') \big) / \A^0 = 0.
\ee
Hence $\widetilde{\kappa}_0(\widetilde{\theta}_l(X'),-))$ is zero, for all $X'$. On the other hand, $\widetilde{\kappa}_0$ is a section of $\wedge^2 \A_1 \otimes \A_{-2}$, hence $\widetilde{\theta}_l (\widetilde{\kappa}_0)$ is zero as well. Together this shows that
\be
\partial^* (\widetilde{\theta}_l \bullet \widetilde{\kappa}_0) = 0.
\ee
Now assume that $\wt{\theta}_l$ is a section of $\A_2 \otimes \A_{l-2}$, and let $X'$ be a section of $\A_{-2}$. Then by equations (\ref{alpha}) and (\ref{beta}),
\be
\partial^* (\widetilde{\theta}_l \bullet \widetilde{\kappa}_0)(X') &=& \alpha (\widetilde{\theta}_l \bullet \widetilde{\kappa}_0)(X') + \beta (\widetilde{\theta}_l \bullet \widetilde{\kappa}_0)(X').
\ee
Then for $\{Z^j\}$ a local, constant homogeneity, frame for $\A^1$, and $\{Z_j\}$ the dual frame for $\A / \A^0$,
\be
\alpha (\widetilde{\theta}_l \bullet \widetilde{\kappa}_0)(X') &=& \sum_j \mc{K}\big( Z^j, \widetilde{\theta}_l (\widetilde{\kappa}_0(X',Z_j)) - \widetilde{\kappa}_0(\widetilde{\theta}_l(X'),Z_j) - \widetilde{\kappa}_0(X', \widetilde{\theta}_l(Z_j) \big) \\
&=& - \sum_j \mc{K}\big( Z^j, \widetilde{\kappa}_0(\widetilde{\theta}_l(X'),Z_j) \big) \\
&=& - \alpha (\widetilde{\kappa}_0)(\widetilde{\theta}_l(X')).
\ee
The important point to notice here is that $\alpha$ maps $\widetilde{\kappa}_0$ to $\A_1 \otimes \A_{-1}$ while $\beta$ maps it to $\A_2 \otimes \A_{-2}$. Hence $\partial^* \widetilde{\kappa}_0 = 0$ implies that $\alpha (\widetilde{\kappa}_0) = \beta (\widetilde{\kappa}_0) = 0$. Hence the above expression must vanish.

Similarly,
\be
\beta (\widetilde{\theta}_l \bullet \widetilde{\kappa}_0)(X') &=& -\frac{1}{2} \Big( \sum_j \widetilde{\theta}_l (\widetilde{\kappa}_0(\mc{K}(Z^j,X')_-,Z_j)) \\
&& - \widetilde{\kappa}_0(\widetilde{\theta}_l(\mc{K}(Z^j,X')_-),Z_j) - \widetilde{\kappa}_0(\mc{K}(Z^j,X')_-,  \widetilde{\theta}_l(Z_j)) \Big). \\
\ee
The bottom two terms must vanish, however -- the first one because $\mc{K}(Z^j,X)_-$ is of homogeneity $\geq -1$, so $\widetilde{\theta}_l(\mc{K}(Z^j,X)_-) = 0$, and the second because the $Z_l$ are of constant homogeneity: if $Z_j$ is of homogeneity $-1$, $\widetilde{\theta}_l(Z_j) = 0$, while if $Z_j$ is of homogeneity $-2$, $Z^j$ must be of homogeneity $+2$, and $\mc{K}(Z^j,X)_- = 0$ (since the expression for $\beta$ does not depend on the choice of local frame, these results are true in general). Hence:
\be
\beta (\widetilde{\theta}_l \bullet \widetilde{\kappa}_0)(X') &=& -\frac{1}{2} \sum_j \widetilde{\theta}_l (\widetilde{\kappa}_0(\mc{K}(Z^j,X')_-,Z_j)) \\
&=& \widetilde{\theta}_l \big(\beta(\wt{\kappa}_0)(X') \big) \\
&=& 0.
\ee

To summarise, $\partial^* (\widetilde{\theta}_l \bullet \widetilde{\kappa}_0)$ maps $\A_{2} \otimes \A_{l-2}$ to $\A_{1} \otimes \A_{l-1}$ and maps this second bundle to zero. Hence, splitting $(\A^1 \otimes \A)_l$ into the above two spaces, we get a block decomposition:
\be
\partial^* \circ \wt{\Psi} (\wt{\theta}_l) = \left( \begin{array}{cc} \square \wt{\theta}_l & \partial^* (\widetilde{\theta}_l \bullet \widetilde{\kappa}_0) \\ 0 & \square \wt{\theta}_l \end{array} \right).
\ee
Hence $(\partial^* \circ \wt{\Psi})^2 = \square^2$, which is invertible on the image of $\partial^*$. Since the image of $\partial^* \circ \wt{\Psi}$ is contained in the image of $\partial^*$, and since its square is invertible on that space, $\partial^* \circ \wt{\Psi}$ itself must be invertible.
\end{lproof}
The previous Lemma completes the demonstration of the Proposition, after we recall that $\wt{\Psi}$ in homogeneity $l$, is precisely $\Psi$ composed with the soldering form $s_{\theta}$.
\end{proof}

Using the above Proposition, one can normalize $\onab$, homogeneity by homogeneity. First, assume that there is an $l$ such that the curvature $\kappa$ of $\onab$ has the property that $\partial^* \kappa$ is of homogeneity $l$. This is certainly true for $l =0$, by assumption. Then choose a $\theta$ such that
\be
\theta_l = - \Psi^{-1}(\partial^* \kappa)_l,
\ee
implying that
\be
(\partial^* \kappa_{\theta})_l = (\partial^* \kappa)_l + \Psi(\pi_l(\theta)) = 0.
\ee
Hence $\onab + \theta$ is has curvature that is $\partial^*$-closed in all homogeneities below $l+1$, so by induction, we can construct a normal Tractor connection. Since we have only added terms of strictly positive homogeneity, the Tractor connection is still partially regular.

Then only ambiguity in the choice of the $\theta_l$ at each homogeneity is by elements transverse to the image of $\partial^*$. But by the cohomology assumptions of Theorem \ref{theo:high}, the space $(T^* \otimes \A)_l$ splits (\cite{capslo}) into the image of $\partial^*$ plus the image of $\mf{g}_l$ under $\partial$. But for $Y \in \mf{g}^l$, $\partial(Y)$ is the lowest homogeneity component of the action of $\exp{Y} \in P$ on $\onab$, and hence, up to higher homogeneities, simply corresponds to an automorphism of $\onab$. Consequently the normalisation procedure is unique, up to automorphism.

This theorem has an important corollary; since the automorphism groups of Cartan connections is a Lie group, of dimension less than or equal to that of the group $G$:
\begin{cor}
On a manifold $M$, let $H \subset T$ be a generic holomorphic distribution, $(G,P)$ a $|2|$-graded parabolic pair such that $H$ has the same rank and co-rank as the dimensions of $\mf{g}_{-1}$ and $\mf{g}_{-2}$. Then if $H^1(\mf{g}^1,\mf{g})$ vanishes in non-negative homogeneities, there is a open dense subset $N$ of $M$, such that the automorphism group of $H|_N$ is a Lie group of dimension less than or equal to that of $G$.
\end{cor}
And the usual result holds in the real category:
\begin{cor}
On a manifold $M$, let $H \subset T$ be a real distribution sufficiently close to a regular one, $(G,P)$ a $|2|$-graded parabolic pair such that $H$ has the same rank and co-rank as the dimensions of $\mf{g}_{-1}$ and $\mf{g}_{-2}$. Then if $H^1(\mf{g}^1,\mf{g})$ vanishes in non-negative homogeneities, there is a open dense subset $N$ of $M$, such that the automorphism group of $H|_N$ is a Lie group of dimension less than or equal to that of $G$.
\end{cor}

\section{Harmonic curvature and BGG sequences}
Partially regular $|2|$-graded geometries have many similarities with regular geometries, going beyond the existence of a unique normal Cartan connection. For instance since $\wedge^3 gr(T)^* \otimes gr(\A)$ has no homogeneity zero component, $\partial$ is trivially zero on $\kappa_0$. Consequently, $\kappa_0$ is harmonic ($\partial^*$ and $\partial$ closed), and so obstructions to flatness of the geometry still lie in the second cohomology bundle:
\begin{prop}
The lowest homogeneity component of the curvature $\kappa$ of a normal $|2|$-graded Cartan connection is both $\partial^*$ and $\partial$ closed, making it a section of the second cohomology bundle $H^2(gr(T),gr(\A))$.
\end{prop}

Furthermore we have:
\begin{prop}
If $\onab$ is the unique normal Tractor connection defined by the geometry and $\mc{V}$ is any Tractor bundle, then the first two splitting operators in the BGG resolution exist. See \cite{BGG} and \cite{BGG2} and for more details on BGG resolutions.
\end{prop}
\begin{proof}
The existence of the first BGG splitting operator is a triviality. To get the second, we need to be able to invert the map
\be
\Phi \to \partial^* \circ d^{\onab} \Phi,
\ee
for $\Phi$ a section of $T^* \otimes \mc{V}$ of homogeneity $l$. To invert this, we need only consider the homogeneity $l$ component. Exactly as in Proposition \ref{curvature:adjustment}, this can be calculated to be
\be
\square \Phi + \Phi(\kappa_0).
\ee
Then, as in the proof of Proposition \ref{curvature:adjustment}, this can be seen to be an invertible bundle map.
\end{proof}

Higher BGG operators are trickier to get. However, if we are prepared to further restrict the set of possible Levi brackets to open dense sets, we are sure of their existence. This is because to construct them, we must be able to invert a certain bundle map, i.e.~avoid the zeros of a certain determinant. Since the bundle map is invertible for $\kappa_0 = 0$, and since the determinant is an analytic function of $\kappa_0$, this means that it must be invertible for generic $\kappa_0$.

In fact, since the set of finite dimensional Tractor bundles is countable, we have the weak result that:
\begin{prop}
For a $|2|$-graded geometry $(M, H \subset TM, G/P)$, at any point $x \in M$, there exists a set $\mc{U}'$ in the $E_x$, the set of possible Levi brackets at $x$. This $\mc{U}'$ has the following properties: for any Levi bracket $\mc{L}$ such that $\mc{L}_x \in E_x$, there exists a unique normal Tractor connection $\onab$ \emph{around} $x$, and all its splitting operators exists \emph{at} $x$.

Moreover this set $\mc{U}'$ is dense in $E_x$ and its complement is of measure zero.
\end{prop}
\begin{proof}
This $\mc{U}'$ is constructed from $\mc{U}$ by removing a countable collection of sets, all of measure zero.
\end{proof}

\subsection{Prolongation procedure}

Regular geometries allow a prolongation procedure for their BBG sequences. The idea is that, given a Tractor bundle $\mc{V}$, the Tractor connection $\onab$ is replaced with a new connection $\nabla^{\mc{V}}$ that is more closely fitted to the BBG sequence and its splitting operators.

See \cite{mrh-BGG} for a detailed analysis of invariant prolongations; for the purposes of this paper, a brief introduction will suffice. The operator $\partial^*$ can be defined on any space $\wedge^l \mf{g}^1 \otimes V$, for $V$ a $\mf{g}$-module. The formulas are similar to equations (\ref{alpha}) and (\ref{beta}); the $\beta$ term is identical, and the $\alpha$ term replaces the bracket action of $Z^l$ with the standard action on $V$:
\be
(\alpha \Psi )(X) = \sum_l Z^l \cdot \Psi (Z_l,X),
\ee
and similarly for higher wedge products. Since this $\partial^*$ is also $P$-invariant, it extends to a bundle map on the associated bundles:
\be
\partial^*: \wedge^l TM^* \otimes \mc{V} \to \wedge^{l-1} TM^* \otimes \mc{V}.
\ee
The tractor derivative $\onab$ defines a differential map going in the opposite direction:
\be
d^{\onab}: \Gamma (\wedge^l TM^* \otimes \mc{V}) \to \Gamma(\wedge^{l+1} TM^* \otimes \mc{V}).
\ee
Now define $\mc{Z}_l \subset \wedge^l TM^* \otimes \mc{V}$ as the kernel of $\partial^*$ and $\mc{B}_l$ as the image of $\partial^*$. Since $(\partial^*)^2 = 0$, $\mc{B}_l \subset \mc{Z}$. The cohomology bundles $\mc{H}_l$ can then naturally be defined as $\mc{Z}_l / \mc{B}_l$. The projection $\mc{Z}_l \to \mc{H}_l$ will be designated by $\pi_l$ (note that this terminology is slightly different from \cite{mrh-BGG}, which uses $\mc{H}$, $\mc{Z}$ and $\mc{B}$ to designate the space of sections of the above bundles).

The standard BGG construction (\cite{BGG} and \cite{BGG2}) rests on the existence of a unique splitting operator $L_l: \Gamma (\mc{H}_l) \to \gamma(\mc{Z}_l)$ such that for any section $s$ of $\mc{H}_l$:
\begin{itemize}
\item $\pi_l \circ L_l (s) = s$ and
\item $d^{\onab} \circ L_l(s) \in \mc{Z}_{l+1}$.
\end{itemize}
Though useful, this construction has the drawback that the corresponding diagrams do not commute; i.e.~in general $d^{\onab} \circ L_l \neq L_{l+1} \circ \pi_{l+1} (d^{\onab}\circ L_l)$. If we define $D_l: \mc{H}_l \to \mc{H}_{l+1}$ as the operator $\pi_{l+1} (d^{\onab}\circ L_l)$, then in general
\be
L_{l+1} \circ D_l \neq d^{\onab} \circ L_l.
\ee

The aim of the prolongation procedure was to construct a unique alternative connection $\wt{\nabla} + \Theta$ on $\mc{V}$, with $\Theta \in \Gamma(T^*\otimes\mf{gl}(\mc{V}))$, chosen so that
\begin{itemize}
\item $\Theta$ is of strictly positive homogeneity and $\Theta(s) \in \Gamma(\mc{B}_1)$ for all sections $s$ of $\mc{V}$,
\item if $\wt{L}_l$ are the splitting operators associated with $\wt{\nabla}$, then $L_0 = \wt{L}_0$ and similarly $\wt{D}_0 = D_0$,
\item $\wt{L}_{1} \circ D_0 = d^{\wt{\nabla}} \circ L_0$.
\end{itemize}
In other words, $D_0$, $d^{\wt{\nabla}}$, $L_0$ and $\wt{L}_1$ form a commuting diagram on the relevant space of sections. Note that $\wt{\nabla}$ is generally \emph{not} a Tractor connection.

To extend this result to the partially regular category, some terminology is needed. Let $s$ be any section of $\mc{V}$, and let $\Theta$ be a section of $T^* \otimes \mf{gl}(\mc{V})$, of homogeneity $\geq 1$, such that for any $s$, $\Theta(s)$ is a section of $\mc{B}_i$ for some fixed $i$.

The whole construction can be extended to partially regular, normal, $\onab$, if two conditions hold; letting $s$ be a section of $\mc{V}$:
\begin{enumerate}
\item The map $gr(\Theta(s))$ to $\square gr(\Theta(s)) + \partial^*(gr(\Theta(\kappa_0)(s)))$ is an invertible map from $gr(\mc{B}_i)$ to itself,
\item $\partial^* (\kappa_0 (s)) = 0$.
\end{enumerate}
The first statement allows us to construct the prolongation by inducting over the homogeneity $i$, exactly as in the regular case. It can be proved by the usual methods used in this paper, since
\be
\beta (gr(\Theta(\kappa_0)(s))) = gr(\Theta(\beta \kappa_0)(s)) = 0.
\ee
The second statement implies that the non-regularity does not mess up the construction from the beginning -- though homogeneity zero can be corrected exactly as the higher homogeneities, the connection will no longer have $\partial$ as its lowest homogeneity piece, throwing the rest of the construction into doubt.

However:
\begin{lemm}
If the bundle $\mc{V}$ is has three or less graded components, then $\partial^* (\kappa_0 (s)) = 0$.
\end{lemm}
\begin{lproof}
Let $s$ be a section of $gr(\mc{V})$, taking values in the lowest two homogeneities. Then since $\kappa_0$ is a section of $\wedge^2 T^*_{1} \otimes T_{-2}$, and since the action of the $T_{-2}$ piece on $s$ must vanish, then $\kappa_0(s) = 0$.

So now assume that $s$ is in the highest homogeneity. As usual, $\beta \kappa_0 (s) = 0$ since $\beta$ only notices the $\wedge^2 T^*$ component. Thus
\be
\partial^* (\kappa_0(s))(X) &=& \alpha (\kappa_0(s))(X) \\
&=& \sum_l Z^l \cdot \kappa_0(X,Z_l) \cdot s \\
&=& \sum_l \mc{K}(Z^l,\kappa_0(X,Z_l)) \cdot s + \kappa_0(X,Z_l) \cdot Z^l \cdot s \\
&=& (\alpha \kappa_0)(s)(X) + \sum_l \kappa_0(X,Z_l) \cdot Z^l \cdot s \\
&=& \sum_l \kappa_0(X,Z_l) \cdot Z^l \cdot s,
\ee
since $(\alpha \kappa_0) = 0$. But since $s$ is in the top homogeneity and $Z^l$ is of homogeneity $\geq 1$, $Z^l \cdot s = 0$.

\end{lproof}
But do $|2|$-graded geometries have Tractor bundles that fit the bill? Indeed this is the case; let the standard Tractor bundle be the bundle
\be
\mc{T}\ = \ \mc{P} \times_{P} V \ = \ \mc{G} \times_{G} V
\ee
where $V$ is the standard representation of $G$. Then:
\begin{lemm}
The standard Tractor bundle has three graded components for all of the non-exceptional $|2|$-graded geometries, with the exception of conformal-spin geometry \cite{me2grad2}.
\end{lemm}
\begin{lproof}
The grading of a representation is the same as the grading of its complexification, since the grading element will have the same eigenvalues on the complexification. Hence the list of complex standard Tractor bundles from \cite{me2grad2} gives the grading, establishing the result.
\end{lproof}

Thus the general result is:
\begin{theo}
The prolongation procedure works for the standard Tractor bundle for all non-exceptional $|2|$-graded geometries, apart from conformal-spin geometry \cite{me2grad2}.
\end{theo}

\bibliographystyle{amsalpha}
\bibliography{ref}

\end{document}